\documentclass[11 pt]{amsart}

\usepackage[english]{babel}
\usepackage[latin1]{inputenc}
\usepackage[T1]{fontenc}
\usepackage{graphicx}
\usepackage{amsmath}
\usepackage{amssymb}
\usepackage{latexsym}
\usepackage{amsthm}
\usepackage[all]{xy}
\usepackage{stmaryrd}
\usepackage{subfigure}
\usepackage{MnSymbol}
\usepackage{comment} 

\def\be{\begin{equation}}
\def\ee{\end{equation}}
\def\beq{\begin{eqnarray*}}
\def\eeq{\end{eqnarray*}}
\def\Z{\mathbb{Z}}

\newcommand{\pic}[3]{\parbox[c]{#1cm}{\includegraphics[scale=#2]{#3}}}

\newcommand{\cerchio}{\includegraphics[width = .4 cm]{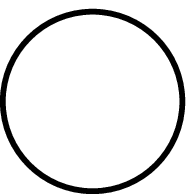}}
\newcommand{\teta}{\includegraphics[width = .4 cm]{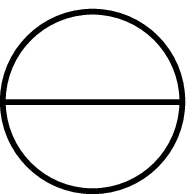}}
\newcommand{\tetra}{\includegraphics[width = .4 cm]{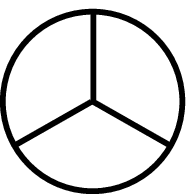}}
\newcommand{\gl}{{\rm gl}}

\newtheorem{theo}{Theorem}[section]
\newtheorem{lem}{Lemma}[section]

\newtheorem{prop}{Proposition}[section]
\theoremstyle{definition}
\newtheorem{defn}{Definition}[section]
\theoremstyle{remark}
\newtheorem{ex}{Example}[section]
\newtheorem{rem}{Remark}[section]

\author{Alessio Carrega}
\address{Dipartimento di Matematica, Largo Pontecorvo 5, 56127 Pisa, Italy}
\email{carrega at mail dot dm dot unipi dot it}

\title[A proof of the shadow formula]{A proof of the shadow formula for the $SU(2)$-Reshetikhin-Turaev-Witten invariant}

\begin{document}

\begin{abstract}
Turaev's \emph{shadow formula} calculates the $SU(2)$-Reshetikhin-Turaev-Witten invariants using shadows, and its expression is somehow similar to a Euler characteristic. We give a short proof of this formula using skein theory.

The formula applies to pairs $(M,G)$ where $M$ is a closed oriented 3-manifold and $G\subset M$ is a (possibly empty) colored framed trivalent graph (for instance, a framed knot or link).

\end{abstract}

\maketitle

\setcounter{tocdepth}{1}
\tableofcontents

\section{Introduction}

In 1984 Vaughan Jones defined his famous polynomial invariant for oriented links using von Newman algebras. This construction revealed many connections between algebra, topology and physics, and it is now a fundamental goal in modern knot theory to ``understand'' the Jones polynomial. After Jones' initial discovery, a variety of new knots and 3-manifolds invariants came out, they are called \emph{quantum invariants} and have many connections with several areas of mathematics and physics.

Some of these invariants arose from representations of braid groups. The images of the generators are examples of \emph{$R$-matrices}, witch play an important role in solving statistical mechanical models and quantum integrable systems in two dimensions. By the end of the 80's, to discover new $R$-matrices Jimbo, Drinfel'd and others developed the formalism of \emph{quantum groups} (or \emph{quantum universal enveloping algebras}), which are deformations $U_q(\mathfrak{g})$ of the universal enveloping algebras $U(\mathfrak{g})$ of the semisimple complex Lie algebras $\mathfrak{g}$ (see for instance \cite{Chari-Pressley, Kassel}). This theory is a part of a more general and categorical theory, where representations of quantum groups form non trivial examples of \emph{ribbon categories}.

The simplest and most studied quantum group is $U_q(\mathfrak{sl}_2)$. Although it is the simplest case, it is rather general and complicated. The fundamental representation of the $SU(2)$-type quantum group yields the $SU(2)$-\emph{Reshetikhin-Turaev-Witten invariants} of 3-manifolds (see \cite{Turaev}). These form the first suggestion of a 3-manifold invariants and was given by Witten \cite{Witten} as a part of his quantum-fields theoretic explanation of the origin of the Jones polynomial. They were then rigorously constructed by Reshetikhin and Turaev \cite{Reshetikhin-Turaev} via surgery presentations of 3-manifolds and quantum groups.

An alternative approach for quantum invariants is provided by \emph{skein theory}. Skeins were introduced by Conway in 1970 for his model of the Alexander polynomial, and Conway's idea became really useful after the work of Kauffman who re-defined the Jones polynomial in a simple combinatorial way. The Jones polynomial is in some sense the simplest quantum invariant, but it is also possible to reproduce all the invariants coming from the representations of $U_q(\mathfrak{sl}_2)$ using skein theory, without referring to quantum groups. This combinatorial description led to many interesting and quite easy computations. The skein formalism was used by Lickorish \cite{Lickorish1, Lickorish2, Lickorish3, Lickorish4}, Blanchet, Habegger, Masbaum and Vogel \cite{BHMV}, and Kauffman and Lins \cite{Kauffman-Lins}, to re-interpret and extend some of the methods above. In particular Lickorish used skein theory to re-define the $SU(2)$-Reshetikhin-Turaev-Witten invariants $I_r(M,G)$ for pairs $(M,G)$ where $M$ is a closed oriented 3-manifolds and $G\subset M$ is a colored framed knotted trivalent graph in $M$. The graph $G$ may for instance be empty or be a framed link or knot.

\emph{Shadows} are 2-dimensional polyhedral objects related to smooth 4-manifolds: these are the 4-dimensional analogue of spines of 3-manifolds. More precisely, they are simple 2-dimensional polyhedra locally-flatly embedded in 4-manifolds, and were defined by Turaev \cite{Turaev:preprint, Turaev} and then considered by various authors, see for instance \cite{Burri, Carrega-Martelli, Costantino1, Costantino2, Costantino-Thurston, Costantino-Thurston:preprint, Goussarov, IK, Martelli, Shu, Thurston, Turaev1}.

In \cite{Turaev} Turaev defines shadows and shows how to get the quantum invariants from them through a formula that works in a general context: for any ribbon category, and hence for any quantum group. In this paper we focus on the Reshetikhin-Turaev-Witten invariants and we reprove this formula using skein theory. The shadow formula is the following:
$$
I_r(M,G) = \kappa^{-\sigma(W_X)} \eta^{\chi(X)} \sum_\xi \frac{\prod_f \cerchio_f^{\chi(f)}A_f \prod_v \tetra_v \prod_{v_\partial} \teta_{v_\partial}}
{\prod_e \teta_e^{\chi(e)} \prod_{e_\partial} \cerchio_{e_\partial}^{\chi(e\partial)} } ,
$$
where $X$ is a shadow (with boundary) of $(M,G)$, $\xi$ runs over all $q$-admissible colorings of $X$ and the other coefficients are standard objects depending on the color and the \emph{gleam} of every vertex, edge and region of $X$ (see Sections~\ref{section:skein_spaces}, \ref{section:shadows} and \ref{section:boundary}). The reader may notice that with this formula the invariant $I_r(M,G)$ is similar to the Euler characteristic of something -- thus increasing the hope of a nice categorification.

We note that this formula also applies to the Jones polynomial of links in $S^3$ and $\#_k(S^2\times S^1)$, and a skein-theoretical proof for that was given in \cite{Carrega-Martelli}.

\subsection*{Structure of the paper}
In Section~\ref{section:skein_spaces} we introduce the skein formalism and recall some important identities, referring to \cite{Lickorish, Kauffman-Lins}.

In Section~\ref{section:surgery} we recall the definition of the $SU(2)$-Reshetikhin-Turaev-Witten invariants $I_r(M)$ of an oriented closed 3-manifold $M$ via a surgery presentation of $M$. Again we can refer to \cite{Lickorish}.

In Section~\ref{section:shadows} we introduce shadows and describe Turaev's shadow formula for $I_r(M)$. We prove this formula then in Section~\ref{section:proof}.

In Section~\ref{section:boundary} we extend everything from $M$ to $(M,G)$, where $G$ is a trivalent framed colored graph in $M$.

\subsection*{Acknowledgments}
The author is warmly grateful to Bruno Martelli for his constant support and encouragement.

\section{Skein theory}\label{section:skein_spaces}

We introduce \emph{skeins}, \emph{colored knotted trivalent graphs}, \emph{Temperley-Lieb algebras}, and some identities that we are going to use in the next sections. The main references are \cite{Lickorish, Kauffman-Lins}.

\subsection{Skein spaces} We start by introducing skein spaces. Let $M$ be an oriented 3-manifold.

\begin{defn}
A \emph{framed link} $L$ in $M$ is a closed 1-submanifold equipped with a \emph{framing}. A framing can be defined as a finite collection of disjoint strips in $M$ whose cores are the components of $L$. Each strip can be an annulus or a M\"obius strip: the latter are usually excluded, but we prefer to include them for technical reasons. If there are no M\"obius strips, we say that $L$ is \emph{orientable}.
\end{defn}

An orientable framed link can be represented with a planar diagram using the \emph{black-board framing}. To recover the framing from the diagram it suffices to orient every component of the link and draw a parallel copy of it on the left (or on the right), thus getting a new link ``parallel'' to the first one.

\begin{figure}[htbp]
$$
\pic{1.8}{0.4}{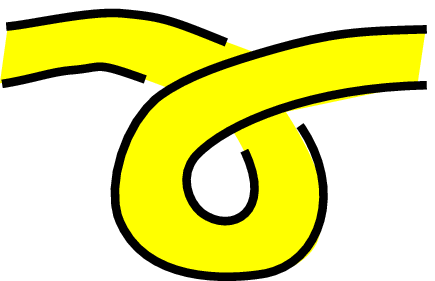} \leftrightarrow \pic{1.8}{0.4}{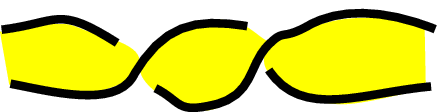} \leftrightarrow \pic{1.8}{0.4}{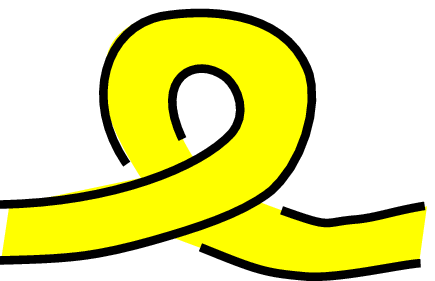} 
$$
\caption{The Reidemeister moves of the first tipe change the framing by adding a full twist.}
\label{figure:framing_change}
\end{figure}

The application of a Reidemester move of the first type changes the framing by adding or removing a full twist (see Fig.~\ref{figure:framing_change}).

We can use the orientation of $M$ to distinguish between a \emph{positive} and a \emph{negative} full twist: the negative twist is the one depicted in Fig.~\ref{figure:framing_change}. Moreover, we can also define a (positive or negative) \emph{half-twist}, which changes the orientability of that component. The composition of two (positive or negative) half-twists is of course a full twist. Two distinct framings on a knot are related by $n$ positive twists, for some $n\in \frac 1 2 \Z$, where if $n<0$ we actually mean $-n$ negative twists.

Each compact surface $S\subset M$ determines an orientable framing on its boundary just by taking a collar of it. Every oriented link $L\subset S^3$ has a \emph{Seifert framing} (or \emph{$0$-framing}), defined by any oriented surface $S$ whose boundary is the oriented $L$. If the link is a knot, the framing does not depend on its orientation.

Now we fix an integer $r\geq 3$ and a primitive $(4r)^{\rm th}$ root of unity $A\in \mathbb{C}$. We mean that $A^{4r}=1$ and $A^n\neq 1$ for each $0< n < 4r$: for instance we might take $A= e^{\pi \sqrt{-1} /2r}$. Although the constant $A$ is often omitted when defining quantum invariants, it is important to note that everything we will say depends on the choice of $A$ and not just on $r$. Furthermore we fix a square root of $A$ and of $-1$ that we denote respectively by $\sqrt A$ and $\sqrt{-1}$ (or $A^{\frac 1 2}$ and $(-1)^{\frac 1 2}$). These two choices almost do not affect the result, it suffices to remember the initial choice and to be coherent.

\begin{defn}\label{defn:skein_space}
Let $M$ be an oriented 3-manifold. Let $V$ be the abstract $\mathbb{C}$-vector space generated by all the framed links in $M$ considered up to isotopy, including the empty set $\varnothing$. The $\mathbb{C}$-\emph{skein vector space} $K_A(M)$ is the quotient of $V$ by the following \emph{skein relations}:
$$
\begin{array}{rcl}
 \pic{1.2}{0.3}{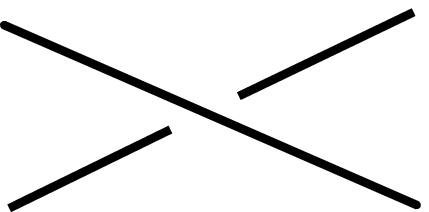}  & = & A \pic{1.2}{0.3}{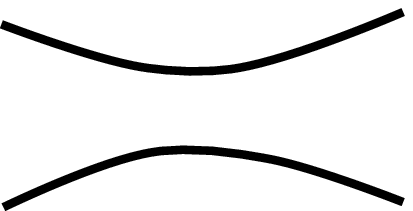}  + A^{-1}  \pic{1.2}{0.3}{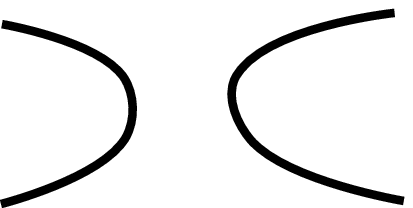}  \\
 D \sqcup \pic{0.8}{0.3}{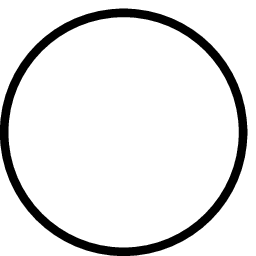}  & = & (-A^2 - A^{-2})  D  \\
K & = & \sqrt{-1} A^{\frac 3 2} K_{-\frac 1 2}
\end{array}
$$
In all relations the links differ only in an oriented 3-ball (we need the orientation of $M$ here). The usual definition includes only the first two relations, and we include here a third one for technical reasons. In the third relation $K$ is any framed knot and $K_{-\frac 12}$ is $K$ with its framing decreased by $-\frac 12$. The relation thus says that making a positive half-twists on any component of a link has the effect of multiplying everything by $\sqrt{-1}A^{\frac 3 2}$.
\end{defn}
The elements of $K_A(M)$ are called \emph{skeins} or \emph{skein elements}.
We can easily deduce that
$$
\pic{0.8}{0.3}{banp.eps} = (-A^2 -A^{-2}) \varnothing.
$$

Kauffman proved that $K_A(S^3)=\mathbb{C}$, generated by the empty link $\varnothing$. Every skein $L$ in $K_A(S^3)$ is equivalent to $L= \langle L\rangle|_A \cdot \varnothing$ where the \emph{Kauffman bracket} $\langle L \rangle$ of $L$ is a Laurent polynomial and $\langle L\rangle|_A$ is its evaluation at the root of unity $A$.

\subsection{Temperley-Lieb algebra}

The $n^{\rm th}$ \emph{Temperley-Lieb algebra} $TL_n$ is a quite famous finitely generated $\mathbb{C}$-algebra. We can think of it as a relative version of the skein space of the 3-cube. The cube has $n$ fixed points on the top side and $n$ on the bottom one and the generators are \emph{framed tangles} from the top to the bottom. Multiplication is defined via superposition of tangles (and extended by linearity). We refer to \cite{Lickorish} for a good definition.

There is a natural injective map $TL_n \rightarrow TL_{n+1}$, obtained by adding a straight line connecting the $(n+1)$'s points. We can identify $TL_n$ with its image.

We also have a natural map $TL_n \rightarrow \mathbb{C}$ called \emph{closure} or \emph{trace}. To get this map we extend by linearity the map defined on tangles that takes an $n$-tangle, identifies the starting points with the corresponding end points, embeds the obtained solid torus in $S^3$ in the standard way, and takes the Kauffman bracket of this skein of $S^3$. 

\begin{figure}[htbp]
$$
1 = \pic{1.5}{0.5}{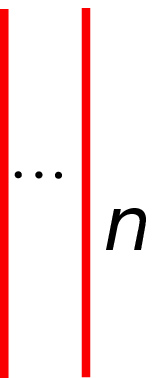} ,\ \ \ e_i = \pic{1.5}{0.5}{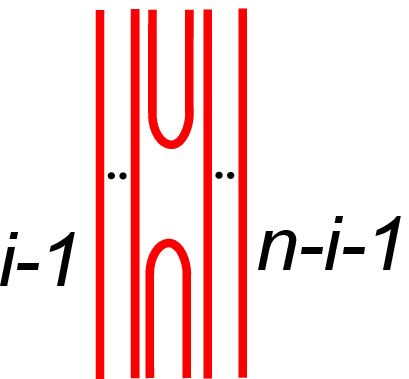}
$$
\caption{Standard generators for the algebra $TL_n$.}
\label{figure:TLgen}
\end{figure}

The Temperley-Lieb algebra $TL_n$ is generated by the elements $1,e_1,\ldots , e_{n-1}$ shown in Fig.~\ref{figure:TLgen}.

\subsection{Colors} We introduce the \emph{Jones-Wenzl projectors} and \emph{colored knotted trivalent graphs}.

\begin{defn}
The \emph{Jones-Wenzl projectors} $f^{(n)}\in TL_n \subset TL_{n+1}$ $0\leq n \leq r-1$, are particular elements of the Temperley-Lieb algebras. Usually the $n^{\rm th}$ projector is denoted as $n$ parallel straight lines covered by a white square box with an $n$ inside. They are completely determined by the following important properties:
\begin{itemize}
\item{$f^{(n)} \cdot e_i = e_i \cdot f^{(n)}=0$ for all $1\leq i \leq n-1$;}
\item{$(f^{(n)} -1)$ belongs on the subalgebra generated by $e_1, \ldots , e_{n-1}$;}
\item{$f^{(n)} \cdot f^{(m)} = f^{(m)}$ for all $n\leq m$.}
\end{itemize}
Moreover the projectors are symmetric, namely they are equivalent to their mirror images under any reflection plane of the cube. This will help below when we color possibly non-orientable framings.
\end{defn}

We denote by $\cerchio_n$ the closure of $f^{(n)}$. It is an element of the field $\mathbb Q(A)$, which is sometimes written as $\Delta_n$ in literature. Since $A$ is a $4r$-primitive root of unity, we have that $\cerchio_n \neq 0$ for $0\leq n \leq r-2$, and $\cerchio_{r-1}=0$. If $A= e^{\frac{\pi\sqrt{-1}}{2r}}$ then
$$
\cerchio_n = \frac{\sin\left( \frac{n\pi}{r}\right) }{\sin\left(\frac{\pi}{r}\right)}.
$$

\begin{defn}
A \emph{framed knotted trivalent graph} $G\subset M$ is a knotted trivalent graph equipped with a \emph{framing}, namely a surface collapsing onto the graph, considered up to isotopy. We admit also closed edges, namely knot components of $G$. A triple $(a,b,c)$ of natural numbers is \emph{admissible}
if they satisfy the triangle inequalities $a\leq b+c$, $b\leq a+c$, $c\leq a+b$, and their sum $a+b+c$ is even. An \emph{admissible coloring} of $G$ is the assignment of a natural number (a \emph{color}) to each edge of $G$ such that the three numbers $a$, $b$, $c$ coloring the three edges incident to each vertex form an admissible triple.

Moreover, a triple of non negative inters $(a,b,c)$ is $q$-\emph{admissible} if it is admissible and $a+b+c\leq 2(r-2)$. In particular it means that each color is at most $r-2$. A coloring of $G$ is $q$-\emph{admissible} if the three numbers $a$, $b$, $c$ coloring the three edges incident to each vertex form a $q$-admissible triple.

A \emph{colored trivalent graph} is a framed knotted trivalent graph with a $q$-admissible coloring.
\end{defn}

\begin{figure}[htbp]
\begin{center}
\includegraphics[width = 7 cm]{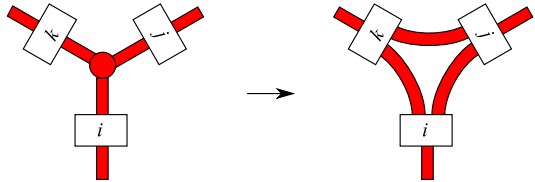}
\end{center}
\caption{A colored framed knotted trivalent graph determines a linear combination of framed links: replace every edge with a Jones-Wenzl projector, and connect them at every vertex via non intersecting strands contained in the depicted bands. For instance there are exactly $\frac{i+j-k}2$ bands connecting the projectors $i$ and $j$.}
\label{figure:trivalent_vertex}
\end{figure}

There is a standard way to define a skein element associated to a colored framed knotted trivalent graph, which agrees with the above definition on framed links with all components colored by $1$. The admissibility requirements on colors allows to associate uniquely to the colored graph a linear combination of framed links by putting the $k^{\rm th}$ Jones-Wenzl projector at each edge colored with $k$ and by substituting vertices with bands as shown in Fig.~\ref{figure:trivalent_vertex}.

\begin{figure}[htbp]
\begin{center}
\includegraphics[width = 9.5 cm]{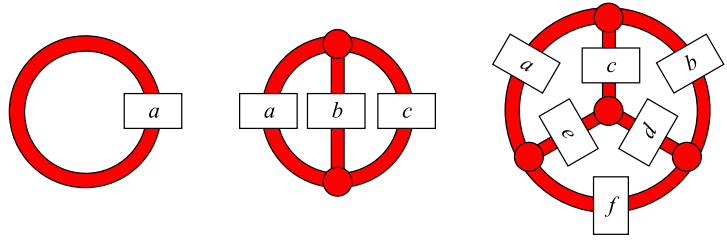}
 \end{center}
\caption{Three important planar colored framed trivalent graphs in $S^3$.}
\label{figure:graphs}
\end{figure}

Three basic planar colored framed trivalent graphs $\cerchio$, $\teta$, and $\tetra$ in $S^3$ are shown in Fig.~\ref{figure:graphs}. Their skeins are some rational functions in $q=A^2$. We can find their expressions in \cite{Kauffman-Lins} or \cite{Carrega-Martelli}.

\subsection{Important objects and identities}

Here we introduce some important objects and identities of the skein spaces that we need. We can find their proofs in \cite{Lickorish}. The first identity is the following

$$
\pic{2}{0.7}{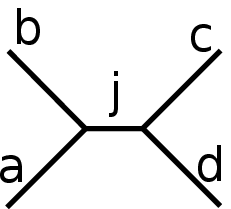} = \sum_i \left\{\begin{matrix} a & b & i \\ c & d & j \end{matrix}\right\} \pic{2}{0.7}{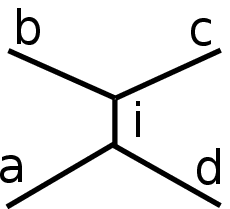}
$$

where the sum is taken over all $i$ such that the colored graphs shown are $q$-admissible. The coefficients between brackets are called $6j$-\emph{symbols}. We show two more identities:
$$
\left\{\begin{matrix} a & b & i \\ c & d & j \end{matrix}\right\} = \frac{ \cerchio_i}{\teta_{a,d,i} \teta_{c,b,i}} \pic{2}{0.7}{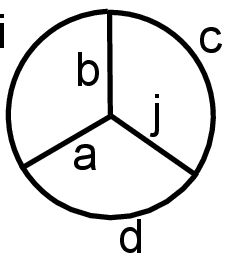}  \quad \quad \quad \quad \text{\cite[Page 155]{Lickorish}},
$$

$$
\pic{1.4}{0.7}{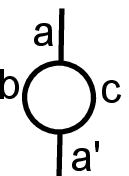} = \left\{\begin{array}{cl}
\frac{\teta_{a,b,c}}{\cerchio_a} \ \pic{1.4}{0.7}{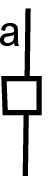} & \text{if } a=a' \\
0 & \text{if } a\neq a'
\end{array}\right. .
$$

\begin{figure}[htbp]
$$
\pic{2}{0.9}{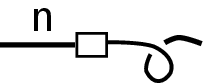} = (-1)^n A^{n^2 + 2n} \ \pic{2}{0.9}{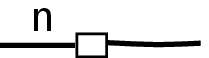}
$$
\caption{A full twist.}
\label{figure:framingchange}
\end{figure}

The effect of a full twist is shown in Fig.~\ref{figure:framingchange}.

More generally, if we change the framing of a framed trivalent graph by adding $k\in \frac 1 2 \Z$ positive twists on an edge colored with $n$, we get the skein of the previous graph multiplied by $(-1)^{nk} A^{nk(n+2)}$:
\begin{prop}\label{prop:half-framingchange}
$$
\pic{2}{0.9}{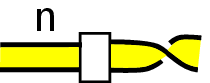} = \sqrt{-1}^n A^{\frac{n^2 + 2n}{2}} \ \pic{2}{0.9}{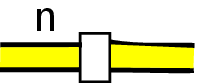}
$$
\begin{proof}
The case for $n=1$ is true for the third skein relation (see Definition~\ref{defn:skein_space}). Then we proceed by induction.
\beq
\pic{2}{0.9}{half-framingchange1.eps} & = & \pic{2.2}{0.9}{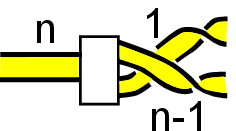} \\
 & = & \sqrt{-1} A^{\frac 3 2} \ \pic{2.2}{0.9}{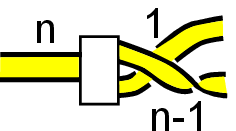} \\
 & = & \sqrt{-1} A^{n+ \frac 1 2} \ \pic{2.2}{0.9}{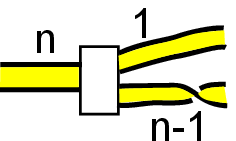} \\
 & = & \sqrt{-1} A^{n+ \frac 1 2} \ \pic{2}{0.9}{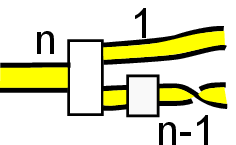} \\
 & = & \sqrt{-1}^n A^{\frac{n^2 + 2n}{2}} \ \pic{2}{0.9}{half-framingchange2.eps}
\eeq
We get the third equality by using $n-1$ times the first skein relation and the fact that the multiplication of a Jones-Wenzl projector with an element of the standard set of generators of $TL_n$ is $0$: $f^{(n)} \cdot e_i = 0$ for $1\leq i \leq n-1$. In fact this property allows us not to consider all the terms multiplied by $A^{-1}$ coming from the application of the first skein relation.
\end{proof}
\end{prop}

\begin{defn}
Consider the skein vector space $K_A(S^1 \times D^2)$ of the solid torus. Let
$\phi_n\in K_A(S^1\times D^2)$ be the core of the solid torus, with trivial framing and colored with $n$. We now construct a particular element of $K_A(S^1\times D^2)$:
$$
\Omega := \eta\sum_{n=0}^{r-2} \cerchio_n \phi_n ,
$$
where
$$
\eta := \left( \sum_{n=0}^{r-2} \cerchio_n^2 \right)^{-\frac{1}{2}} = \frac{A^2-A^{-2}}{\sqrt{-2r}}  \quad \quad \quad \quad \text{\cite[Page 141]{Lickorish}}.
$$
If $A= e^{\frac{\pi\sqrt{-1}}{2r}}$ then
$$
\eta = \sqrt{\frac{2}{r}} \sin\left(\frac{\pi}{r} \right).
$$

If $K$ is a framed knot, we denote by $\Omega K$ the skein obtained by substituting $K$ with $\Omega$. Let $U$, $U_+$, and $U_-$ be the unknot in $S^3$ with framing respectively $0$, $1$, and $-1$. We have

$$
\Omega U = \eta^{-1} \quad \quad \quad \quad \text{\cite[Page 141]{Lickorish}}
$$
and we define
$$
\kappa := \Omega U_+ = \frac{\sum_{n=1}^{4r} A^{n^2}}{2r\sqrt{-2} A^{3+r^2}} \quad \quad \quad \quad \text{\cite[Lemma 14.3]{Lickorish}}.
$$
If $A= e^{\frac{\pi\sqrt{-1}}{2r}}$ then we get
$$
\sum_{n=1}^{4r} A^{n^2} = 2\sqrt{2r} e^{\frac{\pi\sqrt{-1}}{4} } , \ \kappa = \frac{-\sqrt{-1}}{\sqrt{r}} e^{\frac{-\pi\sqrt{-1}}{4r}(2r^2-r +6) }  \quad \quad \text{\cite[Page 148]{Lickorish}}.
$$

Moreover it turns out that
$$
\kappa^{-1} = \Omega U_- \quad \quad \quad \quad \text{\cite[Lemma 13.7]{Lickorish}}.
$$
\end{defn}

\begin{figure}[htbp]
$$
\pic{3}{0.8}{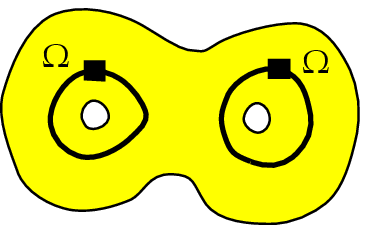} \ = \ \pic{3}{0.8}{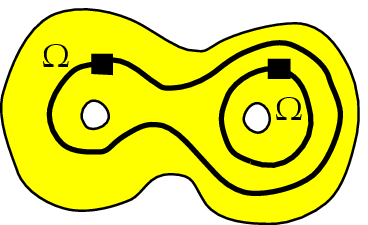}
$$
\caption{Handleslide property for pairs of $\Omega$'s.}
\label{figure:handleslideOmega}
\end{figure}

There is a fundamental relation about $\Omega$ that is called the \emph{handleslide property} for pairs of $\Omega$'s and is shown in Fig.~\ref{figure:handleslideOmega} \cite[Lemma 13.5]{Lickorish}. It says that a Kirby move of the second type does not change the skein, provided that all components are colored with $\Omega$. This move substitutes a component with its bend sum with another component.

\begin{figure}[htbp]
$$
\pic{1.9}{0.8}{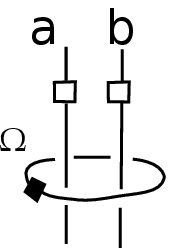} = \left\{\begin{array}{cl}
\frac{\eta^{-1}}{\cerchio_a} \ \pic{1.9}{0.8}{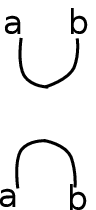} & \text{if }a=b \\
 0 & \text{if }a\neq b
\end{array}\right.
$$
\caption{The 2-strand fusion identity.}
\label{figure:2fusion}
\end{figure}

\begin{figure}[htbp]
$$
\pic{1.9}{0.8}{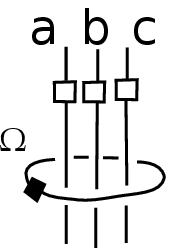} = \left\{\begin{array}{cl}
\frac{\eta^{-1}}{\teta_{a,b,c}} \ \pic{1.9}{0.8}{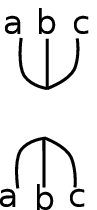} & \text{if }(a,b,c)\text{ is $q$-admissible} \\
 0 & \text{if }(a,b,c)\text{ is not $q$-admissible}
\end{array}\right.
$$
\caption{The 3-strand fusion identity.}
\label{figure:3fusion}
\end{figure}

Two more properties that we will use are the 2- and 3-\emph{strand fusion identity} shown in Fig.~\ref{figure:2fusion} and \ref{figure:3fusion} \cite[Page 159]{Lickorish}. In the left-hand side of the identities we have 2 or 3 strands colored with the $a^{\rm th}$, $b^{\rm th}$ and $c^{\rm th}$ projector. These strands are encircled by a 0-framed unknot colored with $\Omega$.

The 2-strand identity says that if $a=b$ then the left skein is equivalent to $\eta^{-1}\cerchio_a^{-1}$ times the skein obtained by removing the circle, breaking the strands and connecting them in the other way. Otherwise it is equivalent to $0$. The 3-strand identity says that if the triple $(a,b,c)$ is $q$-admissible then the left skein is equivalent to $\eta^{-1}\teta_{a,b,c}^{-1}$ times the skein obtained by removing the circle, breaking the strands and connecting the 3 parts on the same side in a vertex. Otherwise it is equivalent to $0$.

\section{Surgery presentations} \label{section:surgery}

Given an orientable framed link $L$ in a closed oriented 3-manifold $M$ we can construct another closed oriented 3-manifold called the \emph{Dehn surgery} on $L$, in the following way. For each component $L_i$ of $L$ we remove the interior of a closed tubular neighborhood $N_i \cong D^2\times S^1$ such that the framing of $L_i$ corresponds to $\{x\} \times S^1$, and then we glue a solid torus $V_i\cong S^1\times D^2$ to the boundary of $N_i$ via a diffeomorphism of the boundaries that sends a meridian  $ \{y\} \times S^1$ of $V_i$, to the framing of $L_i$.

\begin{defn}
Let $M$ and $N$ be two closed 3-manifolds. A \emph{surgery presentation} of $M$ in $N$ is an orientable framed link $L \subset N$ such that $M$ is obtained from $N$ by Dehn surgery on $L$.
\end{defn}

\begin{rem}
Each surgery presentation in $S^3$ has a 4-dimensional interpretation. In fact with a Dehn surgery we can build not only a 3-manifold $M_L$, but also a 4-manifold $W_L$ whose boundary is $M_L = \partial W_L$. It suffices to see $S^3$ as the boundary of the 4-ball, and then attach a 4-dimensional 2-handle $B_i \cong D^2\times D^2$ along the boundary of a tubular neighborhood $N_i \cong D^2\times S^1$ of each component $L_i$ of the link in the way described above. In fact the boundary of $B_i$ is the union of $N_i$ and $V_i$.
\end{rem}

\begin{theo}[Lickorish, Wallace]
Every orientable closed 3-manifold has a surgery presentation in $S^3$.
\end{theo}

There are two important moves on framed links called  \emph{Kirby moves}. The first one consists in adding a new separated component $U_{\pm}$, that is an unknot with framing $\pm 1$. This corresponds to the connected sum with $S^3$ or (in the 4-dimensional interpretation) to the connected sum with $\mathbb{CP}^2$. The second one is exactly the handleslide that we described above in Fig.~\ref{figure:handleslideOmega}. In the 4-dimensional interpretation, it corresponds to sliding a 2-handle over another. Both Kirby moves do not change the presented 3-manifold.

\begin{theo}[Kirby]
Two surgery presentations $L$ and $L'$ in $S^3$ of the same 3-manifold $M$, are related by isotopies and Kirby moves.
\end{theo}

\begin{defn}
Let $L$ be a surgery presentation in $S^3$ of the closed orientable 3-manifold $M$. The \emph{Reshetikhin-Turaev-Witten} invariant of $M$ (as defined by Lickorish \cite{Lickorish}) is:
$$
I_r(M) := \eta \kappa^{-\sigma(L)} \Omega L,
$$
where $\Omega L$ is the skein element obtained by attaching $\Omega$ to each component of $L$, and $\sigma(L)$ is the signature of the linking matrix of $L$ (the matrix composed by the linking numbers of the components of $L$ and with the framings on the diagonal).
\end{defn}

The quantity $\sigma(L)$ is equal to the signature of the 4-manifold obtained attaching to $D^4$ a 2-handle along each component of $L$. The complex number $I_r(M)$ is a topological invariant. In fact it is clearly invariant by isotopies of links, from the handleslide property we get the invariance under the second Kirby move, and the factor $\kappa^{-\sigma(L)}$ ensures the invariance under the first Kirby move.

\begin{ex}
The connected sum $\#_g(S^1\times S^2)$ of $g$ copies of $S^1\times S^2$ is presented in $S^3$ by the unlink with $g$ components, each one with the 0-framing. Hence $\sigma(L)=0$ and $\Omega L = \eta^{-g}$. Therefore
$$
I_r(S^3)= \eta , \ I_r(S^1\times S^2)= 1 , \ I_r(\#_g(S^1\times S^2)) = \eta^{1-g} .
$$
\end{ex}
\begin{ex}
The unknot with framing $n\in\Z$, $U_n$,  presents the lens space $L(n,1)$. Its linking matrix is the $1\times1$ matrix $(n)$. Hence $\sigma(U_n)={\rm sgn}(n)$. Using equality in Fig.~\ref{figure:framingchange} $n$ times we get
$$
I_r(L(n,1)) = \eta^2 \kappa^{-{\rm sgn}(n)} \sum_{a=0}^{r-2} \cerchio_a^2 (-1)^{an} A^{an(a+2)} .
$$
\end{ex}

\section{Shadows}\label{section:shadows}

We now introduce Turaev's \emph{shadows}, following \cite{Turaev} and \cite{Costantino0}.

\subsection{Generalities}\label{subsec:generalities}

\begin{figure}[htbp]
\begin{center}
\includegraphics[width = 12 cm]{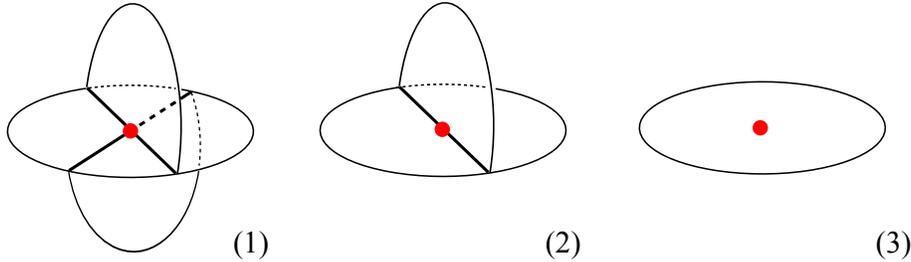}
\caption{Neighborhoods of points in a simple polyhedron.}
\label{figure:models}
\end{center}
\end{figure}

A \emph{simple polyhedron} (without boundary) $X$ is a 2-dimensional compact polyhedron such that each point has a neighborhood homeomorphic to one of the three types (1-3) shown in Fig.~\ref{figure:models}. The three types form subsets of $X$ whose connected components are called \emph{vertices} (1), \emph{edges} (2), and \emph{regions} (3). An edge is either an open segment or a circle; a region is a (possibly non-orientable) connected surface.

\begin{defn}
Let $W$ be a compact oriented 4-manifold with boundary. A \emph{shadow} for $W$ is a simple polyhedron $X\subset W$ such that the following hold:
\begin{itemize}
\item{$X$ is embedded in the interior of $W$;}
\item{$X$ is locally flat, namely every point $p\in X$ has a neighborhood $U$ in $W$ diffeomorphic to $B^3 \times (-1,1)$ and $U\cap X$ is contained in $B^3 \times 0$ as in Fig.~\ref{figure:models};}
\item{$W$ collapses onto $X$.}
\end{itemize}
\end{defn}

\begin{theo} [Turaev]\label{theorem:aaa}
An oriented compact 4-manifold admits a shadow if and only if it admits a handle-decomposition with just 0-, 1- and 2-handles.
\begin{proof}
As we will see in Remark~\ref{rem:W WX} a triangulation of each region of a shadow gives a handle-decomposition of the 4-manifold with just 0-, 1- and 2-handles. On the other hand from such a handle-decomposition we can build a shadow, see \cite{CostantinoPhD}.
\end{proof}
\end{theo}

Let $X$ be a shadow of $W$. Every region $R\subset X$ is equipped with a half-integer called \emph{gleam}, which generalizes the Euler number of closed surfaces embedded in oriented 4-manifolds. The gleam is defined as follows.

The boundary $\partial R$ of $R$ (or more precisely of the closure of $R$) consists of some circles. The shadow $X$ provides an interval sub-bundle of the normal bundle of $\partial R$ in $W$. (The boundary $\partial R$ is actually only immersed in general, but all these definitions work anyway.)

Let $R'$ be a generic small perturbation of $R$ with $\partial R'$ lying in the interval bundle at $\partial R$. The surfaces $R$ and $R'$ intersect only in isolated points, and we count them with signs:
$$
\gl(R) := \frac 1 2 \# (\partial R \cap \partial R') + \#(R\cap R') \in \frac{1}{2} \Z
$$
The half-integer $\gl(R)$ is the \emph{gleam} of $R$ and does not depend on the chosen $R'$. Note that the contribution of $\partial R\cap \partial R'$ on a component of $\partial R$ is even or odd, depending on whether the interval bundle above it is an annulus or a M\"obius strip.

We will say that an abstract simple polyhedron $X$ is an (abstract) \emph{shadow} if each region $R$ of $X$ is equipped with an half-integer $g$, called \emph{gleam}, such that $g \in \Z$ if and only if the interval bundle over $\partial R$ described by $X$ has an even number of non-orientable components.

\begin{theo}[Reconstruction Theorem]\label{theorem:WX}
From an abstract shadow $X$ we can construct a compact oriented 4-manifold $W_X$ such that $X$ is a shadow of $W_X$ and the gleams of $X$, as abstract objects, coincide with the ones given by the embedding in $W_X$.
\begin{proof}
See \cite{Costantino0} or \cite{Turaev}.
\end{proof}
\end{theo}

\begin{rem}\label{rem:W WX}
Let $X$ be a shadow of the oriented compact 4-manifold $W$. Since $W$ collapses onto $X$, the manifold $W$ is diffeomorphic to a regular neighborhood of $X$. This regular neighborhood has a natural decomposition given by the combinatorics of $X$: for each vertex we have a 4-ball, and for each edge we have a $D^3\times  I$ or a $D^3$-bundle over $S^1$, the former attached to the 4-balls as a 1-handle. The 1-skeleton hence thickens to a union of oriented 4-dimensional handlebodies (made of 0- and 1-handles only). Finally, for each region $R$ we have an oriented $D^2$-bundle over $R$, which is $R \times D^2$ when $R$ is orientable, attached to the handlebodies. When $R$ is a disk, this is a 2-handle.
\end{rem}

A shadow $X$ is said to be \emph{standard} if each stratum of $X$ (vertices, edges and regions) is a cell. This implies that each edge is an interval with two different end points and the singular graph is connected.

\begin{rem}\label{rem:standard shadow}
If $X$ is standard the thickening described in the previous remark is a handle-decomposition of $W_X$, where the vertices of $X$ thicken to the 0-handles, the edges to the 1-handles, and the regions to the 2-handles.
\end{rem}

\begin{prop}
Each 4-manifold that admits a shadow has a standard shadow too.
\begin{proof}
The construction in \cite{CostantinoPhD} for Theorem~\ref{theorem:aaa} produces a standard shadow.
\end{proof}
\end{prop}

\subsection{State sum}

Let $X$ be a shadow of $W$. A $q$-\emph{admissible coloring} $\xi$ for $X$ is the assignment of a natural number in $\{0,1,\ldots, r-2\}$ to each region of $X$ (called a \emph{color}), such that for every edge of $X$ the colors of the three incident regions form a $q$-admissible triple.

The \emph{evaluation} of the coloring $\xi$ is the following complex number:
$$
|X|_\xi = \frac{\prod_f \cerchio_f^{\chi(f)}A_f \prod_v \tetra_v }
{\prod_e \teta_e^{\chi(e)} }.
$$
Here the product is taken on all regions $f$, edges $e$, vertices $v$. The symbols
$$
\cerchio_f,\ \teta_e,\ \tetra_v
$$
denote the skein elements of these planar graphs (equipped with a planar framing) in $K_A(S^3)=\mathbb{C}$, colored respectively as $f$ or as the regions incident to $e$ or $v$. For an edge $e$, $\chi(e)$ is the Euler characteristic of its closure in $X$: $\chi(e)=1$ if $e$ is a segment adjacent to two different vertices, otherwise $\chi(e)=0$. In the same way for a region $f$, $\chi(f)$ denotes the Euler characteristic of the closure of $f$ in $X$. We can add to the formula a decorative $\chi(v)$ as exponent of the contribute of a vertex $v$ remembering that $\chi(v)=1$.

The \emph{phase} $A_f$ is the following number:
$$
A_f = (-1)^{gc}A^{-gc(c+2)} ,
$$
where $g$ and $c$ are respectively the gleam and the color of $f$.

Given a shadow $X$ we define
$$
|X|^r := \sum_{\xi} |X|_\xi
$$
where the sum is taken over all the $q$-admissible colorings of $X$. Note that since $X$ is compact and each color of the region is at most $r-2$, the sum is finite.

\subsection{Quantum invariants} Here we state the shadow formula for closed oriented 3-manifolds without graphs (Theorem~\ref{theorem:shadow_formula}). The proof of the formula is postponed to the next section.

We denote with $\sigma(W)$ the signature of a oriented 4-manifold $W$ and with $\chi(X)$ the Euler characteristic of $X$. If $X$ is a shadow of $W$ we have
$$
\chi(X) = \sum_v 1 - \sum_e \chi(e) + \sum_f \chi(f) = \chi(W_X)
$$
where the sum is taken over all the regions $f$, edges $e$ and vertices $v$ of $X$, and $\chi(e)$ and $\chi(R)$ denote respectively the Euler characteristic of the closure in $X$ of the edge $e$ and the region $R$.

\begin{theo}[Shadow formula]\label{theorem:shadow_formula}
Let $X$ be a shadow of the compact oriented 4-manifold $W$. Then
$$
I_r(\partial W) = \kappa^{-\sigma(W)} \eta^{\chi(X)} |X|^r
$$
\end{theo}

\subsection{Bilinear form and signature} Here we show how to get the signature and the intersection form of an oriented 4-manifold from a shadow. 

Let $X$ be a shadow with all regions orientable. Let $R$ be a region of $X$ and let $h\in H_2(X,\Z)$. We denote with $\langle h | R\rangle \in\Z$ the image of $h$ under the map $H_2(X,\Z) \rightarrow H_2(X/(X\setminus R), \Z) \cong \Z$ induced by the quotient identifying the complement of $R$ to a point. The group $H_2(X/(X\setminus R) , \Z)$ is identified with $\Z$ once given an orientation to $R$. The map $Q_X:H_2(X,\Z) \times H_2(X,\Z) \rightarrow \frac 1 2 \Z$ is the bilinear form so defined:
$$
Q_X(h_1,h_2):= \sum_R \langle h_1|
 R \rangle \cdot \langle h_2 | R\rangle  \cdot \gl(R) ,
$$
where the sum runs over all the regions of $X$. It does not depend on the choice of an orientation of the regions. We call $\mathcal{Q}_X$ the \emph{bilinear form} of $X$. We denote with $\sigma(X)$ the signature of $\mathbb{R} \otimes_{\frac 1 2 \Z} Q_X$. We call it the \emph{signature} of $X$.

\begin{theo}\label{theorem:bil_and_sign}
Let $X$ be a shadow of the oriented 4-manifold $W$, and let $f_*:H_2(X,\Z) \rightarrow H_2(W,\Z)$ be the isomorphism induced by the inclusion. Then for any $h_1,h_2\in H_2(X,\Z)$
$$
f_*(h_1) \cdot f_*(h_2) = Q_X(h_1,h_2) ,
$$
where the product on the left-hand side is the intersection product in $H_2(W,\Z)$. Hence 
$$
\sigma(X) = \sigma(W) .
$$
\begin{proof}
See \cite[Section {\rm IX}.5.1.]{Turaev}.
\end{proof}
\end{theo}

\begin{ex}
Let $S$ be the closed orientable surface of genus $g$. Let $M$ be the boundary of the oriented $D^2$-bundle $W$ over $S$ with Euler number $n\in\Z$. If $S=S^2$ $M$ is the lens space $L(n,1)$. The surface $S$ equipped with gleam $n$ is a shadow of $W$. 
$$
\mathcal{Q}_X: \Z \times \Z \rightarrow \Z, \ \ \mathcal{Q}_X(k,h)= k h n . 
$$
By Theorem~\ref{theorem:bil_and_sign} 
$$
\sigma(W) ={\rm sgn}(n) .
$$
Therefore by the shadow formula
$$
I_r(M) = \kappa^{-{\rm sgn}(n)} \eta^{2-2g} \sum_{a=0}^{r-2} \cerchio_a^{2-2g} (-1)^{na} A^{-na(a+2)} .
$$
\end{ex}

\section{Proof of the shadow formula}\label{section:proof} In this section we prove the shadow formula for closed oriented 3-manifolds, namely Theorem~\ref{theorem:shadow_formula}. We will follow the following steps:
 \begin{enumerate}
\item{from a triangulation of the regions of $X$ we construct a surgery presentation $(L,f)$  of $M=\partial W$ in $S^3$, where $f$ is the framing and $L$ is the underlying link;}
\item{we change the framing of the surgery presentation using the identity in Proposition~\ref{prop:half-framingchange};}
\item{we apply the 2-strand fusion identity (Fig.~\ref{figure:2fusion}) to remove the triangulation of the regions;}
\item{we apply the 3-strand fusion identity (Fig.~\ref{figure:3fusion}) to reduce ourselves to a trivalent graph in a tubular neighborhood of a tree (a union of trees);}
\item{we reduce the trivalent graph to isolated tetrahedra;}
\item{we note that the contributions of all our moves give the shadow formula.}
\end{enumerate}

\subsection{Step 1}
We can suppose that $W$, and hence also $X$, is connected. Let $\Gamma \subset S^3$ be an embedding of the 1-skeleton of $X$. Note that we may choose $\Gamma$ so that it has only unknotted edges. 

Select a small 3-ball $B_v$ centered in every vertex $v$ of $\Gamma$. The intersection of one of such 3-balls with $\Gamma$ consists of four unknotted and unlinked strands with one end in $v$ and one in the boundary $\partial B_v$ of the ball. For each vertex $v$ of $X$, put a tetrahedron in $\partial B_v$ whose vertices coincide with $\partial B_v \cap \Gamma$, and give it the framing of the 2-sphere $\partial B_v$. Every non closed edge $e$ of $\Gamma$ connects two different vertices $p_{e,1}$, $p_{e,2}$ of the set of framed tetrahedra (the vertices hold in the same tetrahedron if and only if $\chi(e)=0$). We connect the local framed strands incident to $p_{e,1}$ to the ones of $p_{e,2}$ via three strips running in a regular neighborhood $H_e$ of $e$ as follows. 

\begin{figure}[htbp]
\begin{center}
\includegraphics[width = 9 cm]{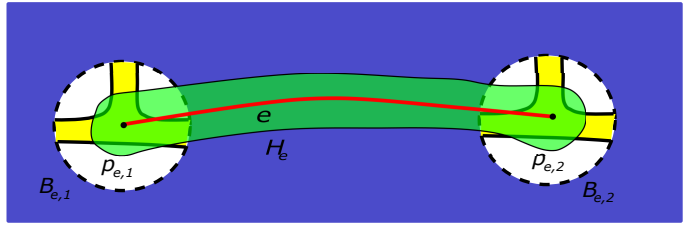}
\caption{The regular neighborhood $H_e$ of a non closed edge $e$ of $\Gamma$. The tube $H_e$ is colored with green, $e$ with red, The framed tetrahedra with yellow, and the blue part covers the rest.}
\label{figure:bridge1}
\end{center}
\end{figure}

\begin{rem}\label{rem:bridge}
Usually a framed graph in $S^3$ with orientable framing, is represented by a diagram in $D^2$ (or $S^2$) giving it the blackboard-framing. Our pictures are made thinking to that method. Note that the regular neighborhoods $H_e$'s are bridges maybe passing over a part of the previous diagram (Fig.~\ref{figure:bridge1}). We will describe the strips using their projection into a rectangle inside $H_e$. This rectangle is exactly the planar one in Fig.~\ref{figure:bridge1}. In this way we can easily construct a diagram of the final framed link. 
\end{rem}

\begin{figure}[htbp]
\begin{center}
\includegraphics[width = 9 cm]{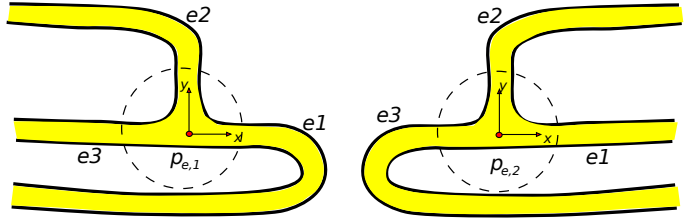}
\caption{Parametrized neighborhood of the vertices.}
\label{figure:attaching_an_edge1}
\end{center}
\end{figure}

The points $p_{e,1}$ and $p_{e,2}$ are trivalent vertices of the set of framed tetrahedra. Let $B_{e,1}$ and $B_{e,2}$ be two small ball neighborhoods of them. We can positively parametrize $B_{e,1}$ and $B_{e,2}$ as in Fig.~\ref{figure:attaching_an_edge1}: the vertex is in the origin of $\mathbb{R}^3$, the framing lies on the plane $\{ (x,y,0) \ | \ x,y\in \mathbb{R} \} $, the first strand lies on the first positive semi-axis ($x\geq 0$, $y=z=0$), the second strand lies on the second positive semi-axis ($y\geq 0$, $x=z=0$), and the third one lies on the first negative semi-axis ($x\leq 0,$ $ y=z=0$).

\begin{figure}[htbp]
\begin{center}
\includegraphics[width = 4.5 cm]{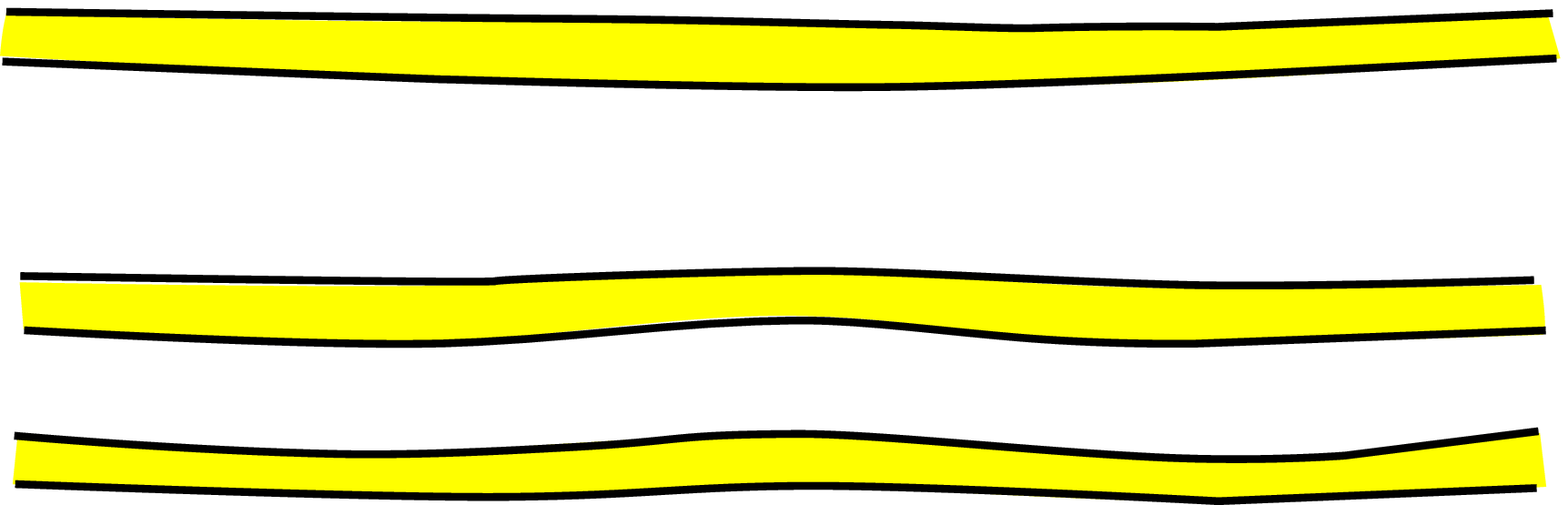} \hspace{1cm} \includegraphics[width = 4.5 cm]{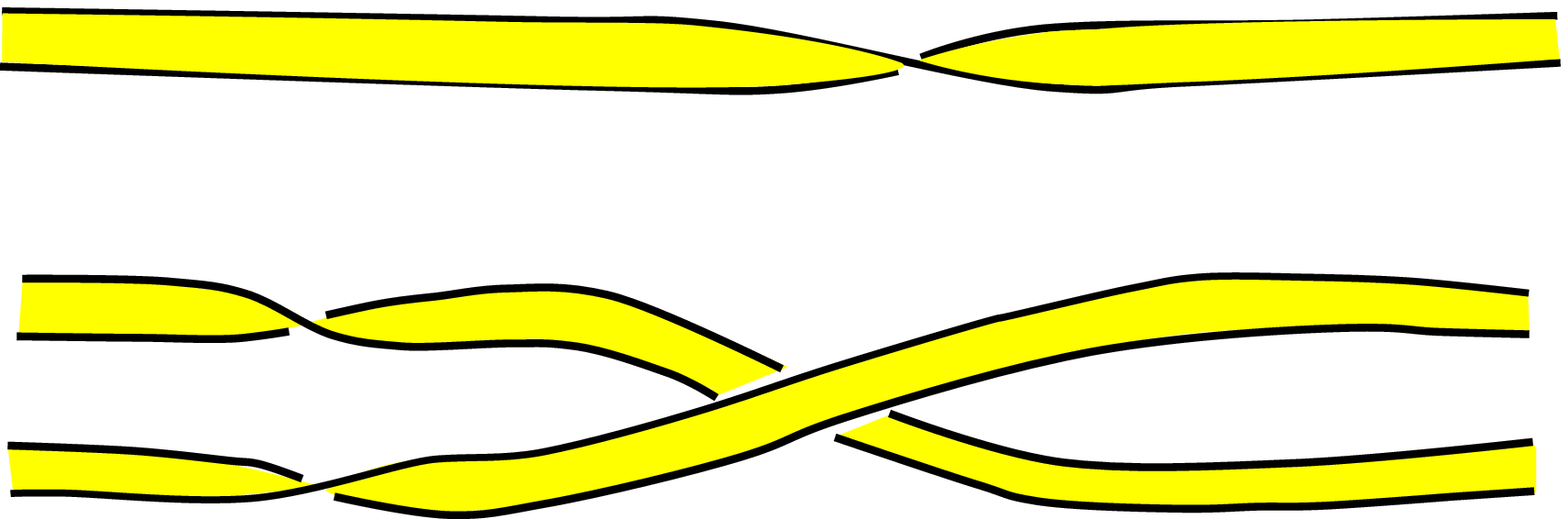}
\caption{Strips connecting the framed strands around the end points of a non closed edge of $\Gamma$.}
\label{figure:attaching_an_edge_b}
\end{center}
\end{figure}

The regular neighborhood of $e$ in $X$ describes a bijection between the framed strands incident to $p_{e,1}$ and the ones incident to $p_{e,2}$.  Up to enumeration of the incident strands, and up to isotopy, we have just two possible bijections:
\begin{itemize}
\item{the second edge is fixed while the first and the third ones are exchanged; }
\item{all the enumerated edges are fixed.}
\end{itemize}
In the first case we connect those strands with three strips passing through $H_e$ and running around $e$ as described in Fig.~\ref{figure:attaching_an_edge_b}-(left): the strips lie in a rectangle whose intersection with $B_{e,1}$ and $B_{e,2}$ is the plane $\{z=0\}$ (see Remark~\ref{rem:bridge}).  In the second case we connect those strands with three strips passing through $H_e$ and running around $e$ as described in Fig.~\ref{figure:attaching_an_edge_b}-(right): the strips can be drawn in that way in a rectangle whose intersection with $B_{e,1}$ and $B_{e,2}$ is the plane $\{z=0\}$ (see Remark~\ref{rem:bridge}). 

Every closed edge $e$ of $\Gamma$ has a solid torus $V_e$ as regular neighborhood.  We divide it in a 3-ball $B_p$ centered in a point $p\in e$ and in a 1-handle $H_e$, $V_e=B_p \cup H_e$. The graph $\Gamma$ intersects $B_p$ in a unknotted properly embedded arc. Put a $\teta$-graph in $\partial B_p$ whose vertices coincides with $\partial B_p \cap \Gamma$, and give it the framing of the 2-sphere $\partial B_p$. Now we have the same situation as before, and we apply the method above to connect the framed strands incident to $\partial B_p \cap \Gamma$ with strips running around the core of $H_e$.

\begin{figure}[htbp]
\begin{center}
\includegraphics[width = 9 cm]{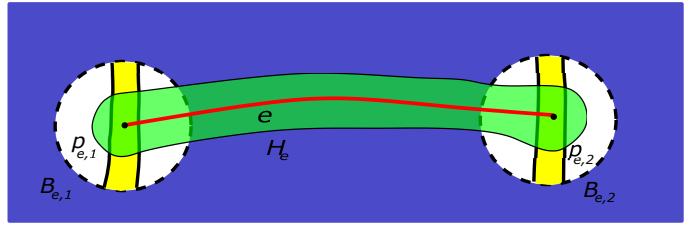}
\caption{The regular neighborhood $H_e$ of an edge $e$ of $\Gamma_R$. The tube $H_e$ is colored with green, $e$ with red, the framed tetrahedra with yellow, and the blue part covers the rest.}
\label{figure:bridge2}
\end{center}
\end{figure}

Now we have obtained a framed link $L'\subset S^3$ (maybe with non orientable framing) lying in a regular neighborhood of $\Gamma$. Let $X^t$ be the simple polyhedron $X$ equipped with the further structure of a triangulation of each region. There is a natural bijectction between the components of $L'$ and the set of connected components of the boundary of the regions of $X$ (every component has the information of an ambient region, thus a component is taken twice if it is in the boundary of two different regions), or, alternatively, with the set of connected components of the boundary of the regular neighborhood of the 1-skeleton of $X$. Let $R$ be a region of $X$. Let $L'_R$ be the sublink of $L'$ whose components are in bijection with the connected components of $R$ ($L' = \cup_R L'_R$, $R_1\neq R_2$ implies $L'_{R_1} \cap L'_{R_2} = \varnothing$). Let $\Gamma_R$ be an embedding of the 1-skeleton of the triangulation of $R$ whose restriction to $\partial R$ is $L'_R$ (without framing). With \emph{internal vertex} of $\Gamma_R$ we mean the image under the embedding of a vertex of the triangulation that lies in the interior of $R$. Select a small 3-ball neighborhood $B_v$ of every internal vertex $v$ of $\Gamma_R$. The set $B_v \cap \Gamma_R$ is a finite number of unknotted and unlinked strands with one end in $v$ and one in $\partial B_v$. Put a 0-framed unknot in $\partial B_v$ containing $\partial B_v \cap \Gamma_R$. Every internal edge $e$ of $\Gamma_R$ (an edge adjacent to an internal vertex) connects two distinct points, $p_{e,1}$ and $p_{e,2}$, lying in framed strands: either in $L'_R$ or in a 0-framed unknot around an internal vertex. Let $H_e$ be an its regular neighborhood (see Fig.~\ref{figure:bridge2}). As before we connect the local framed strands incident to $p_{e,1}$ and $p_{e,2}$ with two strips running through $H_e$. 

\begin{figure}[htbp]
\begin{center}
\includegraphics[width = 9 cm]{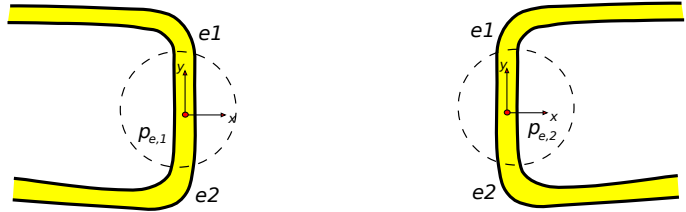}
\caption{Parametrized neighborhood of the non singular points.}
\label{figure:attaching_an_edge4}
\end{center}
\end{figure}

Let $B_{e,1}$ and $B_{e,2}$ be two small ball neighborhoods of $p_{e,1}$ and $p_{e,2}$. We can positively parametrize $B_{e,1}$ and $B_{e,2}$ as in Fig.~\ref{figure:attaching_an_edge4}: the point is in the origin of $\mathbb{R}^3$, the framing lies on the plane $\{ (x,y,0) \ | \ x,y\in \mathbb{R} \} $, the first strand lies on the second positive semi-axis ($y\geq 0$, $x=z=0$) and the second one lines on the second negative semi-axis ($y\leq 0$, $x=z=0$).

\begin{figure}[htbp]
\begin{center}
\includegraphics[width = 4.5 cm]{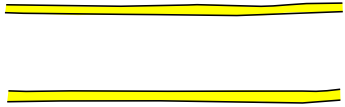}
\end{center}
\caption{Strips connecting the framed strands around the end points of an internal edge of $\Gamma_R$.}
\label{figure:attaching_an_edge_c}
\end{figure}

The regular neighborhood of $e$ in $X$ describes a bijection between the edges incident to $p_{e,1}$ and the ones incident to $p_{e,2}$.  Up to enumeration of the incident edges, and up to isotopy, we have just one possible bijection: the enumerated edges are fixed. Thus we connect these strands with two strips passing through $H_e$ and running around $e$ as described in Fig.~\ref{figure:attaching_an_edge_c}: the strips lie in a rectangle whose intersection with $B_{e,1}$ and $B_{e,2}$ is the plane $\{z=0\}$ (see Remark~\ref{rem:bridge} and Fig.~\ref{figure:bridge2}).

Repeat this procedure for every region $R$ of $X$. Now we have obtained a framed link lying in a regular neighborhood of a connected graph $\Gamma^t \subset S^3$ containing $\Gamma$. The graph $\Gamma^t$ is an embedding of the 1-skeleton of $X^t$ extending the one of $X$ and the ones of the regions, $\Gamma, \Gamma_R \subset \Gamma^t$. If $\Gamma^t$ is a circle, add a 0-framed meridian unknot that encircles the framed link. By ``\emph{maximal tree}'' we mean a collapsible subgraph of an ambient graph whose set of vertices is the same as the ambient one. If $\Gamma^t$ has a vertex (is not a circle), take a maximal tree $T$ of $\Gamma^t$. For each edge of $\Gamma^t$ not lying in $T$ add a 0-framed unknot that encircles the strips running around its regular neighborhood.

We get a framed link $(L,f_1)$ in $S^3$ (the framing $f_1$ may not be orientable) with some components, $\epsilon_1,\ldots, \epsilon_k$, corresponding to the boundary components of the tubular neighborhood of the 1-skeleton of $X^t$, and the other ones, $\delta_1,\ldots,\delta_g$, are the added 0-framed unknots. For each region $R$ of $X$ take the curves $\epsilon_{i(R,1)}, \ldots , \epsilon_{i(R,k_R)}$ that lie in the regular neighborhood of $\Gamma_R$. Add to their framings some half-twists so that the sum (with sign) of the added twists is equal to $-\gl(R)\in \frac 1 2 \Z$. So we get the framing $f$.


\begin{rem}
The framing $f$ clearly depends on the choice made when we modified $f_1$. If $X$ is standard we can take the trivial triangulation: just one triangle for each region. In that case there is only one choice to modify $f_1$ and we can easily see that it makes $f$ orientable:

The gleam$\gl(R)$ is an integer if and only if there is an even number of edges of $X$ that form $\partial R$ and define the second bijection in the list given when we described the connecting strips. Only in that case we applied a half-twist to each strip (see Fig.~\ref{figure:attaching_an_edge_b}-(right)). Hence $L'_R$ is a knot with an orientable framing plus an even number of half-twists. Namely $f$ is orientable.  
\end{rem}

\begin{theo}
One of the choices above makes $(L,f)$ a surgery presentation of $M$ in $S^3$.
\end{theo}

\begin{rem}
Links in $\#_g(S^1\times S^2)$ like $\epsilon_1 \cup \ldots \cup \epsilon_k$, are called \emph{universal links}. The word ``universal'' is due precisely to the fact that we can get any orientable connected closed 3-manifold by surgering on them. See \cite{Costantino-Thurston} and \cite{Costantino2}, in particular \cite[Proposition 3.35]{Costantino-Thurston} and \cite[Proposition 3.36]{Costantino-Thurston}.
\end{rem}

\subsection{Step 2}
Let $(L,f)$ be the surgery presentation in $S^3$ of $\partial W$ described above. We remind that
$$
I_r(\partial W) = \eta \kappa^{-\sigma((L,f))} \Omega (L,f) ,
$$
where $\Omega (L,f)$ is the skein element got by coloring each component of $(L,f)$ with $\Omega= \eta \sum_{n=0}^{r-2}\cerchio_n \phi_n$. The skein $\Omega (L,f)$ is equal to the sum over $0 \leq n_1, \ldots n_k \leq r-2$ of $\eta^k \prod_i \cerchio_{n_i}$ times the skein element obtained by giving to every $\delta_j$ the color $\Omega$ and to each $\epsilon_i$ the $n_i^{\rm th}$ projector. Fix one of these colorings.

We change again the framing of the $\epsilon_{i(R)}$'s by adding $\gl(R)$ positive twists, namely we return to $f_1$. For Proposition~\ref{prop:half-framingchange} each of these framing changes produces a multiplication by $A_R$ ($ = (-1)^{\gl(R)n_{i(R)}} A^{-\gl(R) n_{i(R)}(n_{i(R)}+2)}$), the phase of the colored region $R$.

\begin{rem}
$f$ depends on a made choice and not all gives a surgery presentation. However we can easily check with this use of Proposition~\ref{prop:half-framingchange} that the skein elements given by all these $f$'s are the same.
\end{rem}

\subsection{Step 3}
Near each $\delta_j$ encircling two strands, we have the situation of the 2-strand fusion (see Fig.~\ref{figure:2fusion}). Apply the identity to each such $\delta_j$. The identity splits $L$ to a new link, and multiplies it by a coefficient. Furthermore the identity says that the colors of the $\epsilon_i$'s related to the same region must be equal, otherwise the summand is null. The new colored framed link consists of the curves $\delta_{j,1},\ldots , \delta_{j,g'}$ colored with $\Omega$ and the colored curves $\epsilon'_1,\ldots , \epsilon'_{k'}$ encircled by the $\delta_{j,l}$'s, all with the framing induced by $f_1$. Applying an isotopy we can see that the framed link $\epsilon'_1 \cup \ldots \cup \epsilon'_{k'}$ is equal to the framed link $L'$ we constructed before introducing the triangulations of the regions.

\subsection{Step 4}
Near each $\delta_{j,l}$ we have the situation of the 3-strand fusion (see Fig.~\ref{figure:3fusion}). Hence we apply it for each $l=1,\ldots,g'$. Thus if there is a non $q$-admissible triple $(n_{i1}, n_{i2}, n_{i3})$ of colors of $\epsilon'_i$'s encircled by a $\delta_{j,l}$, then the summand is null. After the application of all the 3-strand fusions it remains an unknotted 0-framed trivalent graph $G$ in the regular neighborhood a tree. 

\subsection{Step 5}
$G$ has two vertices for each edge of $X^t$ which does not lie in $T$ and three parallel edges for each edge of $T$. A tree graph has vertices connected by an edge with only another vertex, the \emph{leaves}, and the other ones are connected with two different vertices. Near a leaf of $T$, $G$ has the form of the left-hand side of the equality in the following Lemma~\ref{lem:a}. We apply the equality to the parts of $G$ near a fixed leaf of $T$. We get a multiple of another trivalent graph $G_1$ that is made in the same way of $G$ but encircling the embedding of the subgraph $T_1$ of $T$ obtained removing the fixed leaf and the adjacent edge. We repeat this procedure until we finish the edges of $T$. The lemma says also that if there is an edge of $T$ with the three strands along it that are colored with a non $q$-admissible triple, then the summand is null. Therefore our summation is taken over all the $q$-admissible colorings of $X$.

\subsection{Step 6}
Fix one of these $q$-admissible colorings $\xi$. By the  applications of the 2-strand fusion identities and the framing change, we get for each region $R$ a contribute of $A_R \cdot \cerchio_{n_i}^{\chi(R)-1}
$ where $A_R$ is the phase of $R$ colored with $n_i$ (we added also the contribution $\cerchio_{n_i}$ times the considered skein element). 

By the applications of Lemma~\ref{lem:a} we get a contribute $\tetra_v$ for each colored vertex of $X$, and a $\teta_e$ for each colored edge of $X$ lying on the maximal tree of the 1-skeleton of $X^t$.

By the applications of the 3-strand fusion we get a contribute of $\teta_e$ for each colored edge of $X$ that does not lie on the maximal tree.

$\sigma((L,f))=\sigma(W_1)$, where $W_1$ is the 4-manifold obtained attaching to $D^4$ a 2-handle along each component of $(L,f)$. Let $W$ be the 4-manifold obtained giving a dot to the $\delta_j$'s. Let $W_g$ be the 4-dimensional orientable handlebody of genus $g$ (the compact 4-manifold with a handle-decomposition with just $k$ 0-handles and $k-g+1$ 1-handles for some $k>0$), and let $\#_{\partial,g}(D^2\times S^2)$ be the boundary connected sum of $g$ copies of $D^2\times S^2$. They have the same boundary: $\#_g(S^1\times S^2)$. The framed link $(L,f)$ has a corresponding framed link $L'$ in $\#_g(S^1\times S^2)$. Let $W'$ be the 4-manifold obtained attaching to $\#_g(S^1\times S^2)\times [-1,1]$ a 2-handle along each component of $L'\times\{1\}$. We have $W= W'\cup W_g$ and $W_1= W' \cup \#_{\partial,g}(D^2\times S^2)$. The signature is additive and $\sigma(W_g)=\sigma(\#_{\partial,g}(D^2\times S^2)) =0$. Therefore $\sigma((L,f))=\sigma(W_1)= \sigma(W') + \sigma(W_g) = \sigma(W') + \sigma(\#_{\partial,g}(D^2\times S^2)) = \sigma(W)$. 

Moreover the 2- and 3-strand fusions give a contribute of $\eta^{-g}$. Therefore we get that the summand of $\Omega (L,f)$ is equivalent to
$$
\eta^{k-g}|X|_\xi  .
$$
$k$ is the number of triagles of $X^t$ and $g$ is the genus of its 1-skeleton. Hence $1-g+k=\chi(X)=\chi(W)$. Therefore
$$
I_r(\partial W) = \eta^{\chi(W)} \kappa^{-\sigma(L)} |X|^r .
$$

\begin{lem}\label{lem:a}
$$
\pic{1.7}{0.7}{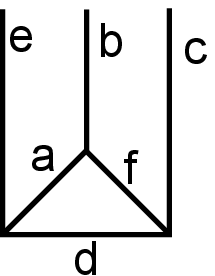} = \left\{\begin{array}{cl}
\frac{\pic{1.7}{0.5}{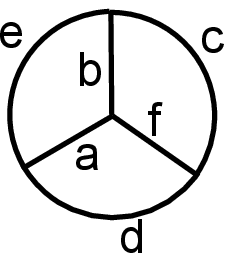}}{\teta_{e,b,c}} \ \ \pic{1.7}{0.7}{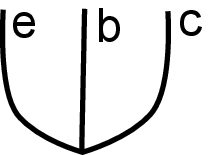} & \text{if }(e,b,c)\text{ is $q$-admissible}\\
0 & \text{if }(e,b,c)\text{ is not $q$-admissible}
\end{array}\right.
$$
\begin{proof}
\beq
\pic{1.7}{0.7}{lemma1.eps} & = & \sum_i \left\{\begin{matrix} a & b & i \\ c & d & f \end{matrix}\right\} \ \pic{1.7}{0.7}{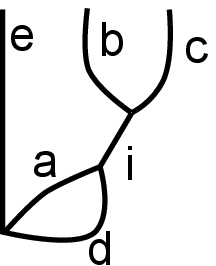} \\
 & = & \sum_{i,j} \left\{\begin{matrix} a & b & i \\ c & d & f \end{matrix}\right\} \left\{\begin{matrix} d & e & j \\ i & d & a \end{matrix}\right\} \ \pic{1.7}{0.7}{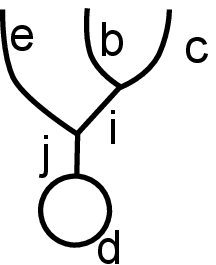}
\eeq
\beq
 & = & \left\{\begin{array}{cl}
 \left\{\begin{matrix} a & b & e \\ c & d & f \end{matrix}\right\} \left\{\begin{matrix} d & e & 0 \\ e & d & a \end{matrix}\right\} \ \pic{1.7}{0.7}{lemma4.eps} & \text{if }(e,b,c)\text{ is $q$-admissible} \\
 0 & \text{if }(e,b,c)\text{ is not $q$-admissible}
 \end{array}\right.\\
 & = & \left\{\begin{array}{cl}
 \frac{\pic{1.7}{0.5}{tetra_color.eps} }{\teta_{e,b,c}} & \text{if }(e,b,c)\text{ is $q$-admissible} \ \pic{1.7}{0.7}{lemma4.eps} \\
 0 & \text{if }(e,b,c)\text{ is not $q$-admissible}
 \end{array}\right.
\eeq
\end{proof}
\end{lem}

\section{Invariant for links and colored trivalent graphs}\label{section:boundary}
In this section we describe the invariant for pairs $(M,G)$ where $M$ is an oriented closed 3-manifold and $G$ an embedded colored trivalent graph in $M$, we introduce the notion of \emph{shadow with boundary} and we give and prove the \emph{shadow formula} (Theorem~\ref{theorem:shadow_formula2}) for this case with boundary.

\subsection{Skein theory}
\begin{defn}
Let $G$ be a $q$-admissible colored framed trivalent graph in the closed oriented 3-manifold $M$, and let $L$ be a surgery presentation of $M$ in $S^3$. Let $G'$ be the colored framed trivalent graph of $S^3$ that corresponds to $G$ by the surgery on $L$. If $r$ is bigger equal than the biggest color of the edges of $G$, we define
$$
I_r(M,G) := \eta \kappa^{-\sigma(L)} \Omega (L,G') ,
$$
where $\Omega(L,G')$ is the skein element obtained by coloring the components of $L$ with $\Omega$ and the edges of $G'$ with the projectors corresponding to the colors.
\end{defn}

\begin{rem}
It is clearly an invariant for $q$-admissible colored framed trivalent graphs. If $G$ has no vertices it is a colored framed link. We get a family of invariants of non colored framed links just by varying $r$ and coloring the link components with $r-2$. If the surgery link $L$ is empty $M=S^3$ we can copy the Kauffman's construction of the Jones polynomial to obtain also a family of invariants for oriented links. It suffices to multiply our invariants for framed links in $S^3$ by $((-1)^{r-2}A^{(r-2)^2+2(r-2)})^{-w(L')}$, where $w(L')$ is the sum of the signs of a diagrammatic representation of $L'$ $(= G')$ (the \emph{writhe number}). In this case, $L=\varnothing$, $M=S^3$, we get an evaluation of the Kauffman bracket or the colored Jones polynomial for each $r$.
\end{rem}

\subsection{Shadows with boundary}
We enlarge the notion of shadow.

\begin{figure}[htbp]
\begin{center}
\includegraphics[width = 12 cm]{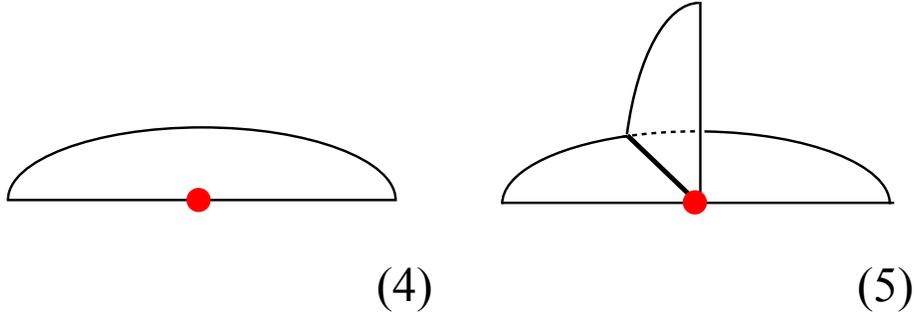}
\caption{Neighborhoods of boundary points in a simple polyhedron.}
\label{figure:models45}
\end{center}
\end{figure}

A \emph{simple polyhedron with boundary} $X$ is a 2-dimensional compact polyhedron such that each point has a neighborhood homeomorphic either to one of the three types (1-3) shown in Fig.~\ref{figure:models}, or to the two types (4-5) shown in Fig.~\ref{figure:models45}. The two types form subsets of $X$ whose connected components are called \emph{external edges} (4) and \emph{external vertices} (5). An external edge is either an open segment or a circle. Together they form a trivalent graph called the \emph{boundary} $\partial X$ of $X$. Regions touching $\partial X$ are called \emph{external regions}.

\begin{defn}
Let $G$ be a colored trivalent graph of the closed oriented 3-manifold $M$. A \emph{shadow} for $(M,G)$ is a simple polyhedron $X\subset W$ such that $W$ is a 4-manifold bounded by $M$ and the following hold:
\begin{itemize}
\item{$X$ is properly embedded in $W$;}
\item{$X$ is locally flat, namely every point $p\in X$ has a neighborhood $U$ in $W$ either diffeomorphic to $B^3 \times (-1,1)$ and $U\cap X$ is contained in $B^3 \times 0$ as in Fig.~\ref{figure:models}, or $U \cong B^2 \times (-1,0] \times (-1,1)$ and $U \cap X$ is contained in $B^2\times (-1,0]\times\{0\}$ as in Fig.~\ref{figure:models45};}
\item{$W$ collapses onto $X$;}
\item{$\partial X= G$.}
\end{itemize}
\end{defn}

We can give the analogous notions and results of Section~\ref{subsec:generalities} for this case with boundary. In particular here we give the notion of \emph{bilinear form} and \emph{signature} of a shadow.
 \subsection{Bilinear form and signature}
Let $X$ be a shadow (maybe with boundary) with only orientable regions. For any $h\in H_2(X,\partial X;\Z)$ $\langle h|R \rangle \in \Z$ is the image of $h$ under the map $H_2(X,\partial X;\Z) \rightarrow H_2(X/(X\setminus R);\Z) \cong \Z$ given by the homomorphism $X/ \partial X \rightarrow X/(X\setminus R)$. To define it we need to fix an orientation of $R$. The map $\tilde{\mathcal{Q}}_X : H_2(X,\partial X;\Z)\times H_2(X,\partial X;\Z) \rightarrow \frac 1 2 \Z$ is the bilinear form so defined: $\tilde{\mathcal{Q}}_X(h_1,h_2)= \sum_R \langle h_1| R\rangle \cdot \langle h_2| R\rangle \cdot \gl(R)$. Here $R$ runs over all the regions and we do not need to fix an its orientation. The group $H_2(X;\Z)$ is contained in $H_2(X,\partial X;\Z)$ and we call $\mathcal{Q}_X$ the restriction of $\tilde{\mathcal{Q}}_X$ to $H_2(X,\Z)$. This is the \emph{bilinear form} of $X$. The \emph{signature} of $X$, $\sigma(X)$, is defined as the signature of $\mathbb{R} \otimes_{\frac 1 2 \Z} \mathcal{Q}_X$. Theorem~\ref{theorem:bil_and_sign} works also for this case with boundary (always with $\mathcal{Q}_X$).

\subsection{Quantum invariants}
Let $X$ be a shadow with colored boundary (in the sense of trivalent graphs). A $q$-\emph{admissible coloring} $\xi$ for $X$ that \emph{extends} the coloring of $\partial X$, is the assignment of a color to each region of $X$ (an integer), such that for every interior edge of $X$ the colors of the three incident regions form a $q$-admissible triple, the color of an external region coincide with the color of its boundary component in $\partial X$, and for every external vertex the triple of colors of the incident external edges in it form a $q$-admissible triple (the coloring of $\partial X$ is $q$-admissible).

The \emph{evaluation} of the coloring $\xi$ is the following function:
$$
|X|_\xi = \frac{\prod_f \cerchio_f^{\chi(f)}A_f \prod_v \tetra_v \prod_{v\partial} \teta_{v\partial}}
{\prod_e \teta_e^{\chi(e)} \prod_{e\partial} \cerchio_{e\partial}^{\chi(e\partial)} } .
$$
Here the product is taken on all regions $f$, inner edges $e$, inner vertices $v$, external edges $e\partial$, external vertices $v\partial$. The symbols
$$
\cerchio_f,\ \cerchio_{e \partial},\ \teta_e,\ \teta_{v\partial}, \  \tetra_v
$$
denote the skein element of these graphs in $K_A(S^3)=\mathbb{C}$, colored respectively as $f$, $e\partial$, or as the regions incident to $e$, $v\partial$, or $v$. The quantity $A_f$ is the phase with the color and the gleam of the region $f$. The quantities $\chi(e)$, $\chi(f)$ and $\chi(e\partial)$ are the Euler characteristic of the closure in $X$ of the inner edge $e$, the region $f$ and the external edge $e\partial$. While $\chi(v\partial)=1$ for every external vertex of $v\partial$.

Let $X$ be a shadow with colored boundary.
$$
|X|^r := \sum_{\xi} |X|_\xi
$$
where the sum is taken over all the $q$-admissible colorings of $X$ that extend the coloring of $\partial X$.

\begin{theo}[Shadow formula]\label{theorem:shadow_formula2}
Let $G$ be a colored framed trivalent graph of the closed oriented 3-manifold $M$, and let $X$ be a shadow of $(M,G)$. If $r$ is bigger equal to the biggest color of $G$ then
$$
I_r(M,G) = \kappa^{-\sigma(W_X)} \eta^{\chi(X)} |X|^r
$$
\begin{proof}
We proceed in the same way of Section~\ref{section:proof}:
\begin{enumerate}
\item{from a triangulation of the regions of $X$ we construct a surgery presentation of $M$ in $S^3$ (a union of copies of $S^3$) together with the colored trivalent graph $G'$ corresponding to $G$;}
\item{we change the framing of the surgery presentation by using the equality in Proposition~\ref{prop:half-framingchange};}
\item{we apply the 2-strand fusion identity Fig.~\ref{figure:2fusion} to remove the triangulation of the regions;}
\item{we apply the 3-strand fusion identity Fig.~\ref{figure:2fusion} to reduce ourselves to a trivalent graph in a tubular neighborhood of a tree (a union of trees);}
\item{we apply Lemma~\ref{lem:a} to reduce the trivalent graph to isolated tetrahedra;}
\item{we note that the contributions of all our moves give the shadow formula.}
\end{enumerate}
Here we give some clarifications in order to do the first step. Take as $\Gamma \subset S^3$ an embedding of the 1-skeleton of $X$ minus the edges adjacent to (or contained in) the boundary $\partial X$. After the procedure for the edges of $\Gamma$, not all the vertices of the isolated tetrahedra are going to be removed. In fact the ones corresponding to internal edges adjacent to $\partial X$ remain. This happens if and only if $G\subset \#_g(S^1\times S^2)$ has vetices. There is still a bijection between the edges (closed or not) of the resulting framed graph $L'\subset S^3$ and a set consisting of all the boundary components of the internal regions and all the components of $G$. Fix a triangulation $X^t$ for each internal region of $X$. As before we use the embeddings $\Gamma_R$'s, but this time $R$ runs just over all the internal regions. The graph $\Gamma^t$ is an embedding of the 1-skeleton of $X^t$ minus the edges (of $X$) adjacent to $\partial X$. The graph $T$ is still a maximal tree of such connected graph $\Gamma^t$. The framed graph $G'\subset S^3$ corresponding to $G\subset \#_g(S^1\times S^2)$ is the subgraph of $L'$ whose components correspond to the components of $G$.

To end we note that in the shadow formula the contribute of the external regions together with the one of the external edges has no effects, and the same happens for the contribute of the edges adjacent to external vertices together with the external vertices.
\end{proof}
\end{theo}

\end{document}